\documentclass[12pt,reqno]{amsart}
\usepackage{}
\usepackage{amsfonts}
\usepackage{a4wide}
\allowdisplaybreaks
\numberwithin{equation}{section}

\newtheorem{theorem}{Theorem}[section]
\newtheorem{proposition}[theorem]{Proposition}

\newtheorem{lemma}[theorem]{Lemma}

\theoremstyle{definition}

\begin{document}

\title
 [fractional Schr\"{o}dinger equations with a general nonlinearity]
 {Ground state solution of fractional Schr\"{o}dinger equations with a general nonlinearity*}\footnotetext{*This work is supported by Natural Science Foundation of China (Grant No. 11601530, 11371159).}

\maketitle
 \begin{center}
\author{Yi He}
\footnote{Corresponding Author: Yi He. Email addresses: heyi19870113@163.com (Y. He).}
\end{center}

\begin{center}
\address{School of Mathematics and Statistics, South-Central University For Nationalities, Wuhan, 430074, P. R. China}
\end{center}

\maketitle

\begin{abstract}
In this paper, we study the following fractional Schr\"{o}dinger equation:
\[
\left\{ \begin{gathered}
  {( - \Delta )^s}u + mu = f(u){\text{ in }}{\mathbb{R}^N}, \hfill \\
  u \in {H^s}({\mathbb{R}^N}),{\text{ }}u > 0{\text{ on }}{\mathbb{R}^N}, \hfill \\
\end{gathered}  \right.
\]
where $m>0$, $N>2s$, ${( - \Delta )^s}$, $s \in (0,1)$ is the fractional Laplacian. Using minimax arguments, we obtain a positive ground state solution under general conditions on $f$ which we believe to be almost optimal.


\end{abstract}

\maketitle

\section{Introduction and Main Result}

\setcounter{equation}{0}

We cunsider the following fractional Schr\"{o}dinger equation:
\begin{equation}\label{1.1}
\left\{ \begin{gathered}
  {( - \Delta )^s}u + mu = f(u){\text{ in }}{\mathbb{R}^N}, \hfill \\
  u \in {H^s}({\mathbb{R}^N}),{\text{ }}u > 0{\text{ on }}{\mathbb{R}^N}, \hfill \\
\end{gathered}  \right.
\end{equation}
where $m>0$, $N>2s$, ${( - \Delta )^s}$, $s \in (0,1)$ is the fractional Laplacian. The nonlinearity $f:\mathbb{R} \to \mathbb{R}$ is a continuous function. Since we are looking for positive solutions, we assume that $f(t)=0$ for $t<0$. Furthermore, we need the following conditions:\\
$(f_1)$ $\mathop {\lim }\limits_{t \to {0^ + }} f(t)/t = 0$;\\
$(f_2)$ $\mathop {\lim }\limits_{t \to  + \infty } f(t)/{t^{2_s^ *  - 1}} = 1$ where $2_s^ *  = 2N/(N - 2s)$;\\
$(f_3)$ $\exists \lambda > 0$ and $2 < q < {2_s^ *}$ such that $f(t) \ge \lambda {t^{q - 1}} + {t^{2_s^ *  - 1}}$ for $t \ge 0$.\\
Note that, for the case $s=1$, $(f_1)$-$(f_3)$ were first introduced by J. Zhang, Z. Chen and W. Zou \cite{zcz}. This hypothesis can be regarded as an extension of the celebrated Berestycki-Lions' type nonlinearity (see \cite{bl1,bl2}) to the fractional Schr\"{o}dinger equations with critical growth.

Equation \eqref{1.1} has been derived as models of many physical phenomena, such as phase transition, conservation laws, especially in fractional quantum mechanics, etc., \cite{fqt}. \eqref{1.1} was introduced by N. Laskin \cite{l2,l3} as an extension of the classical nonlinear Schr\"{o}dinger equations $s=1$ in which the Brownian motion of the quantum paths is replaced by a L\'{e}vy flight. We refer to \cite{dpv} for more physical backgrounds.

In recent years, the study of fractional Schr\"{o}dinger equations has attracted much attention from many mathematicians. In \cite{crs,css,s1}, L. Caffarelli, L. Silvestre $et~al$ investigated free boundary problems of fractional Schr\"{o}dinger equations and obtained some regularity estimates. In \cite{cs1,cs2}, X. Cabr\'{e} and Y. Sire studied the existence, uniqueness, symmetry, regularity, maximum principle and qualitative properties of solutions to the fractional Schr\"{o}dinger equations in the whole space. For more results, we refer to \cite{ap,bcps1,bcps,cw,dpv,fqt,jlx,rs}.

Our main result is as follows:

\begin{theorem}\label{1.1.}
Assume that the nonlinearity $f$ satisfies $(f_1)$-$(f_3)$. If $N \ge 4s$, $2<q<{2_s^ * }$ or $2s< N < 4s$, $4s/(N - 2s) < q < 2_s^ * $, then for every $\lambda > 0$, \eqref{1.1} possesses a positive ground state solution. Moreover, the same conclusion holds provided that $2s< N < 4s$, $2 < q \le 4s/(N - 2s)$ and $\lambda > 0$ sufficiently large.
\end{theorem}

We note that, to the best of our knowledge, there is no result on the existence of positive ground state solutions for fractional Schr\"{o}dinger equation under $(f_1)$-$(f_3)$.

The proof of Theorem~\ref{1.1.} is based on variational method. The main difficulties lie in two aspects: (i) The facts that the nonlinearity $f(u)$ does not satisfy $({\text{AR}})$ condition and the function $f(s)/s$ is not increasing for $s > 0$  prevent us from obtaining a bounded Palais-Smale sequence ((PS) sequence in short) and using the Nehari manifold respectively. (ii) The unboundedness of the domain $\mathbb{R}^N$ and the nonlinearity $f(u)$ with critical growth lead to the lack of compactness.

To complete this section, we sketch our proof.

To treat the nonlocal problem \eqref{1.1}, we use the L. Caffarelli and L. Silvestre extension method \cite{cs} to study a corresponding extension problem
\begin{equation}\label{1.5}
\left\{ \begin{gathered}
   - {\text{div}}({y^{1 - 2s}}\nabla w) = 0{\text{ in }}\mathbb{R}_ + ^{N + 1}, \hfill \\
   - {k_s}\mathop {\lim }\limits_{y \to {0^ + }} {y^{1 - 2s}}\frac{{\partial w}}
{{\partial y}}(x,y) =  - mw + f(w){\text{ on }}{\mathbb{R}^N} \times \{ 0\} . \hfill \\
\end{gathered}  \right.
\end{equation}
with the corresponding functional
\[
{I_m}(w) = \frac{{{k_s}}}
{2}\int_{\mathbb{R}_ + ^{N + 1}} {{y^{1 - 2s}}|\nabla w{|^2}} dxdy + \frac{m}
{2}\int_{{\mathbb{R}^N}} {{w^2}(x,0)} dx - \int_{{\mathbb{R}^N}} {F(w(x,0))} dx,{\text{ }}w \in {X^{1,s}}(\mathbb{R}_ + ^{N + 1}).
\]
where $F(s): = \int_0^s {f(t)} dt$ and ${X^{1,s}}(\mathbb{R}_ + ^{N + 1})$ is defined as the completion of $C_0^\infty (\overline {\mathbb{R}_ + ^{N + 1}} )$ under the norm
\[
{\| w \|_{{X^{1,s}}(\mathbb{R}_ + ^{N + 1})}} = {\Bigl( {\int_{\mathbb{R}_ + ^{N + 1}} {{y^{1 - 2s}}|\nabla w{|^2}} dxdy + \int_{{\mathbb{R}^N}} {{w^2}(x,0)} dx} \Bigr)^{1/2}}.
\]

Motivated by J. Hirata, N. Ikoma and K. Tanaka \cite{hit}, by applying the General Minimax principle (Theorem~2.8 of \cite{w1}) to the composite functional
\[
{I_m} \circ \Phi (\theta ,w): = {I_m}(w({e^{ - \theta }}x,{e^{ - \theta }}y)),{\text{ }}(\theta ,w) \in \mathbb{R} \times {X^{1,s}}(\mathbb{R}_ + ^{N + 1}),
\]
we construct a bounded ${{\text{(PS)}}_{{c_m}}}$ sequence $\{ {w_n}\} _{n = 1}^\infty  \subset {X^{1,s}}(\mathbb{R}_ + ^{N + 1})$ with an extra property ${P_m}({w_n}) \to 0$ as $n \to \infty $ where $c_m$ is the mountain pass level of $I_m$ and ${P_m}(w)=0$ is the Pohozaev's identity of \eqref{1.5} (Proposition~\ref{3.2.} below). Proceeding by standard arguments, the existence of ground state solutions for \eqref{1.5} follows.

This paper is organized as follows, in Section 2, we give some preliminary results. In Section 3, we prove the main result Theorem~\ref{1.1.}.\\

%\medskip

\section{Preliminaries}

\setcounter{equation}{0}
In this section, we collect some preliminary results. Recall that for $s \in (0,1)$,  ${D^s}({\mathbb{R}^N})$ is defined by the completion of $C_0^\infty ({\mathbb{R}^N})$ with respect to the Gagliardo norm
\[
{\| u \|_{{D^s}({\mathbb{R}^N})}} = {\left( {\int_{{\mathbb{R}^{2N}}} {\frac{{|u(x) - u(y){|^2}}}
{{|x - y{|^{N + 2s}}}}} dxdy} \right)^{1/2}}
\]
and the embedding ${D^s}({\mathbb{R}^N}) \hookrightarrow {L^{2_s^ * }}({\mathbb{R}^N})$ is continuous, that is
\[
{\| u \|_{{L^{2_s^ * }}({\mathbb{R}^N})}} \le C(N,s){\| u \|_{{D^s}({\mathbb{R}^N})}}
\]
by Theorem~1 of \cite{ms}. The fractional Sobolev space ${H^s}({\mathbb{R}^N})$ is defined by
\[
{H^s}({\mathbb{R}^N}) = \left\{ {u \in {L^2}({\mathbb{R}^N}):\int_{{\mathbb{R}^{2N}}} {\frac{{|u(x) - u(y){|^2}}}
{{|x - y{|^{N + 2s}}}}dxdy} } \right\}
\]
endowed with the norm
\[
{\| u \|_{{H^s}({\mathbb{R}^N})}} = {\| u \|_{{D^s}({\mathbb{R}^N})}} + {\| u \|_{{L^2}({\mathbb{R}^N})}}.
\]
For $N>2s$, we see from Lemma~2.1 of \cite{ap} that
\begin{equation}\label{c4}
{H^s}({\mathbb{R}^N}){\text{ is continuously embedded in }}{L^p}({\mathbb{R}^N}){\text{ for }}p \in [2,2_s^*].
\end{equation}

An important feature of the operator ${( - \Delta )^s} (0<s<1)$ is its nonlocal character. A common approach to deal with this problem was proposed by L. Caffarelli and L. Silvestre \cite{cs}, allowing to transform \eqref{1.1} into a local problem via the Dirichlet-Neumann map in the domain $\mathbb{R}_ + ^{N + 1}: = \{ (x,t) \in {\mathbb{R}^{N + 1}}:t > 0\} $. For $u \in {D^s}({\mathbb{R}^N})$, the solution $w \in {X^s}(\mathbb{R}_ + ^{N + 1})$ of
\[\left\{ \begin{gathered}
   - {\text{div}}({y^{1 - 2s}}\nabla w) = 0{\text{ in }}\mathbb{R}_ + ^{N + 1}, \hfill \\
  w = u{\text{ on }}{\mathbb{R}^N} \times \{ 0\}  \hfill \\
\end{gathered}  \right.\]
is called $s$-harmonic extension of $u$, denoted by $w = {E_s}(u)$. The $s$-harmonic extension and the fractional Laplacian have explicit expressions in terms of the Poisson and the Riesz kernels, respectively
\begin{equation}\label{2.1}
w(x,y) = P_y^s * u(x) = \int_{{\mathbb{R}^N}} {P_y^s(x - \xi)u(\xi )} d\xi ,
\end{equation}
where
\[
P_y^s(x): = c(N,s)\frac{{{y^{2s}}}}
{{{{(|x{|^2} + {y^2})}^{(N + 2s)/2}}}}
\]
with a constant $c(N,s)$ such that $\int_{{\mathbb{R}^N}} {P_1^s(x)} dx = 1$ (see \cite{jlx}).

Here, the space ${X^s}(\mathbb{R}_ + ^{N + 1})$ is defined as the completion of $C_0^\infty (\overline {\mathbb{R}_ + ^{N + 1}} )$ under the norm
\[
{\| w \|_{{X^s}(\mathbb{R}_ + ^{N + 1})}}: = {\Bigl( {\int_{\mathbb{R}_ + ^{N + 1}} {{k_s}{y^{1 - 2s}}|\nabla w{|^2}} dxdy} \Bigr)^{1/2}}.
\]
From \cite{bcps}, the map ${E_s}( \cdot )$ is an isometry between ${D^s}({\mathbb{R}^N})$ and ${X^s}(\mathbb{R}_ + ^{N + 1})$, i.e. for $w = {E_s}(u)$,
\begin{equation}\label{2.3}
{\| u \|_{{D^s}({\mathbb{R}^N})}} = {\| w \|_{{X^s}(\mathbb{R}_ + ^{N + 1})}}.
\end{equation}
On the other hand, for a function $w \in {X^s}(\mathbb{R}_ + ^{N + 1})$, we shall denote its trace on ${\mathbb{R}^N} \times \{ 0\} $ as $u(x) : = {\text{Tr}}(w) = w(x,0)$. This trace operator is also well defined and it satisfies
\begin{equation}\label{2.4}
{\| u \|_{{D^s}({\mathbb{R}^N})}} \le {\| w \|_{{X^s}(\mathbb{R}_ + ^{N + 1})}}.
\end{equation}

\begin{lemma}\label{2.2.}
(Theorem~2.1 of \cite{bcps}) For every $w \in {X^s}(\mathbb{R}_ + ^{N + 1})$, it holds that
\[
S(s,N){\Bigl( {\int_{{\mathbb{R}^N}} {|u{|^{2_s^ * }}} dx} \Bigr)^{2/2_s^ * }} \le \int_{\mathbb{R}_ + ^{N + 1}} {{y^{1 - 2s}}|\nabla w{|^2}} dxdy,
\]
where $u = {\text{Tr}}(w)$. The best constant takes the exact value
\[
S(s,N) = \frac{{2{\pi ^s}\Gamma (1 - s)\Gamma ((N + 2s)/2)\Gamma {{(N/2)}^{2s/N}}}}
{{\Gamma (s)\Gamma ((N - 2s)/2)\Gamma {{(N)}^{2s/N}}}}
\]
and it is achieved when $u_{\delta}$ takes the form
\[
u_{\delta}(x) = {\delta ^{(N - 2s)/2}}{(|x {|^2} + {\delta ^2})^{ - (N - 2s)/2}}
\]
for some $\delta  > 0$ and $w_{\delta} = {E_s}(u_{\delta})$.
\end{lemma}

\section{Proof of the main results}
In view of \cite{cs}, \eqref{1.1} can be transformed into
\begin{equation}\label{3.1}
\left\{ \begin{gathered}
 - {\text{div}}({y^{1 - 2s}}\nabla w) = 0{\text{ in }}\mathbb{R}_ + ^{N + 1}, \hfill \\
   - {k_s}\mathop {\lim }\limits_{y \to {0^ + }} {y^{1 - 2s}}\frac{{\partial w}}
{{\partial y}}(x,y) =  - mw + f(w){\text{ on }}{\mathbb{R}^N} \times \{ 0\}  \hfill \\
\end{gathered}  \right.
\end{equation}
with the corresponding functional
\[
{I_m}(w) = \frac{{{k_s}}}
{2}\int_{\mathbb{R}_ + ^{N + 1}} {{y^{1 - 2s}}|\nabla w{|^2}} dxdy + \frac{m}
{2}\int_{{\mathbb{R}^N}} {{w^2}(x,0)} dx - \int_{{\mathbb{R}^N}} {F(w(x,0))} dx,{\text{ }}w \in {X^{1,s}}(\mathbb{R}_ + ^{N + 1}).
\]
In view of \cite{cw,rs}, if $w \in {X^{1,s}}(\mathbb{R}_ + ^{N + 1})$ is a weak solution to \eqref{3.1}, the following Pohozaev's identity holds:
\begin{equation}\label{3.2}
{P_m}(w) = \frac{{{k_s}(N - 2s)}}
{2}\int_{\mathbb{R}_ + ^{N + 1}} {{y^{1 - 2s}}|\nabla w{|^2}} dxdy + \frac{{mN}}
{2}\int_{{\mathbb{R}^N}} {{w^2}(x,0)} dx - N\int_{{\mathbb{R}^N}} {F(w(x,0))} dx = 0.
\end{equation}

\begin{lemma}\label{3.1.}
${I_m}$ possesses the Mountain-Pass geometry (see \cite{ar}), i.e.\\
$(i)$ There exist ${\rho _0},{\alpha _0} > 0$ such that ${I_m}(w) \ge {\alpha _0}$ for all $w \in {X^{1,s}}(\mathbb{R}_ + ^{N + 1})$ with ${\| w \|_{{X^{1,s}}(\mathbb{R}_ + ^{N + 1})}} = {\rho _0}$.\\
$(ii)$ $\exists {w _0} \in {X^{1,s}}(\mathbb{R}_ + ^{N + 1})$ such that ${I_m}({w _0}) < 0$.
\end{lemma}

\begin{proof}
$(i)$ By $(f_1)$ and $(f_2)$, $\forall \delta  > 0$, $\exists {C_\delta } > 0$ such that
\begin{equation}\label{3.3}
f(w) \le \delta |w| + {C_\delta }|w{|^{2_s^ *  - 1}}{\text{ and }}F(w) \le \delta |w{|^2} + {C_\delta }|w{|^{2_s^ * }}.
\end{equation}
Choosing $\delta  = m/4$ in \eqref{3.3}, we see from Lemma~\ref{2.2.} that
\[
{I_m}(w) \ge \frac{1}
{4}\| w \|_{{X^{1,s}}(\mathbb{R}_ + ^{N + 1})}^2 - C\| w \|_{{X^{1,s}}(\mathbb{R}_ + ^{N + 1})}^{2_s^ * },
\]
then taking ${\rho _0},{\alpha _0} > 0$ small, $(i)$ holds.

$(ii)$ For $R>0$, $T>0$, we define
\[{w_{R,T}}(x,y) = \left\{ \begin{gathered}
  T,{\text{ if }}(x,y) \in  {B_R^ + (0)} , \hfill \\
  T(R + 1 - {{\text{(}}|x{|^2} + {y^2})^{1/2}}),{\text{ if }}(x,y) \in B_{R + 1}^ + (0)\backslash B_R^ + (0), \hfill \\
  0,{\text{ if }} (x,y) \in \mathbb{R}_ + ^{N + 1}\backslash  {B_{R + 1}^ + (0)} , \hfill \\
\end{gathered}  \right.\]
then ${w_R} \in {X^{1,s}}(\mathbb{R}_ + ^{N + 1})$. By $(f_3)$ and the polar coordinate transformation, we have
\[\begin{gathered}
  {\text{    }}{I_m}({w_{R,T}}(x/\theta ,y/\theta )) \hfill \\
   = \frac{{{k_s}}}
{2}{\theta ^{N - 2s}}\int_{\mathbb{R}_ + ^{N + 1}} {{y^{1 - 2s}}|\nabla {w_{R,T}}{|^2}} dxdy  + {\theta ^N}\Bigl[ {\frac{m}
{2}\int_{{\mathbb{R}^N}} {w_{R,T}^2(x,0)} dx - \int_{{\mathbb{R}^N}} {F({w_{R,T}}(x,0))} dx} \Bigr] \hfill \\
   \le \frac{{{k_s}}}
{2}{\theta ^{N - 2s}}{T^2}\int_{B_{R + 1}^ + (0)\backslash B_R^ + (0)} {{y^{1 - 2s}}} dxdy  + {\theta ^N}\Bigl[ {\frac{m}
{2}\int_{\Gamma _{R + 1}^0(0)} {w_{R,T}^2(x,0)} dx - \frac{1}
{{2_s^ * }}\int_{\Gamma _R^0(0)} {w_{R,T}^{2_s^ * }(x,0)} dx} \Bigr] \hfill \\
   \le C{\theta ^{N - 2s}}{T^2}\int_R^{R + 1} {{r^{N + 1 - 2s}}} dr + {\theta ^N}\Bigl[ {C\Bigl( {\frac{m}
{2}{T^2} - \frac{1}
{{2_s^ * }}{T^{2_s^ * }}} \Bigr){R^N} + C{T^2}({{(R + 1)}^N} - {R^N})} \Bigr] \hfill \\
   \le C{T^2}{R^{N + 1 - 2s}}{\theta ^{N - 2s}} + \Bigl( {C\Bigl( {\frac{m}
{2}{T^2} - \frac{1}
{{2_s^ * }}{T^{2_s^ * }}} \Bigr){R^N} + C{T^2}{R^{N - 1}}} \Bigr){\theta ^N}. \hfill \\
\end{gathered} \]

Choosing a large ${T_0} > 0$ such that $\frac{m}
{2}T_0^2 - \frac{1}
{{2_s^ * }}T_0^{2_s^ * } < 0$, then we can choose a large ${R_0} > 0$ such that $C\Bigl( {\frac{m}
{2}T_0^2 - \frac{1}
{{2_s^ * }}T_0^{2_s^ * }} \Bigr)R_0^N + CT_0^2R_0^{N - 1} < 0$, at last, we select a large $\bar \theta  > 0$ to ensure that ${I_m}({w_{{R_0},{T_0}}}(x/\bar \theta ,y/\bar \theta )) < 0$, ${w_{{R_0},{T_0}}}$ is the desired $w_0$.
\end{proof}

Hence we define the Mountain-Pass level of ${I_m}$:
\begin{equation}\label{3.4}
{c_m}: = \mathop {\inf }\limits_{\gamma  \in {\Gamma _m}} \mathop {\sup }\limits_{t \in [0,1]} {I_m}(\gamma (t)),
\end{equation}
where the set of paths is defined as
\begin{equation}\label{3.5}
{\Gamma_m}: = \left\{ {\gamma  \in C([0,1],{X^{1,s}}(\mathbb{R}_ + ^{N + 1})):\gamma (0) = 0{\text{ and }}{I_m}(\gamma (1)) < 0} \right\}.
\end{equation}
By Lemma~\ref{3.1.}(i), we see that ${c_m} > 0$. Moreover, we denote
\[
{b_m}: = \inf \{ {I_m}(w):w \in {X^{1,s}}(\mathbb{R}_ + ^{N + 1})\backslash \{ 0\} {\text{ be a nontrivial solution of \eqref{3.1}}}\} .
\]

Next, we will construct a (PS) sequence $\{ {w_n}\} _{n = 1}^\infty $ for $I_m$ at the level $c_m$ that satisfies ${P_m}({w_n}) \to 0$ as $n \to \infty $, i.e.
\begin{proposition}\label{3.2.}
There exists a sequence $\{ {w_n}\} _{n = 1}^\infty $ in ${X^{1,s}}(\mathbb{R}_ + ^{N + 1})$ such that, as $n \to \infty $,
\begin{equation}\label{3.6}
{I_m}({w_n}) \to {c_m},{\text{ }}{I'_m}({w_n}) \to 0,{\text{ }}{P_m}({w_n}) \to 0.
\end{equation}
\end{proposition}

\begin{proof}
Define the map $\Phi :\mathbb{R} \times {X^{1,s}}(\mathbb{R}_ + ^{N + 1}) \to {X^{1,s}}(\mathbb{R}_ + ^{N + 1})$ for $\theta  \in \mathbb{R}$, $w \in {X^{1,s}}(\mathbb{R}_ + ^{N + 1})$ and $(x,y) \in \mathbb{R}_ + ^{N + 1}$ by $\Phi (\theta ,w) = w({e^{ - \theta }}x,{e^{ - \theta }}y)$. For every $\theta  \in \mathbb{R}$, $w \in {X^{1,s}}(\mathbb{R}_ + ^{N + 1})$, the functional ${I_m} \circ \Phi $ is computed as
\[\begin{gathered}
  {I_m} \circ \Phi (\theta ,w) = \frac{{{k_s}}}
{2}{e^{(N - 2s)\theta }}\int_{\mathbb{R}_ + ^{N + 1}} {{y^{1 - 2s}}|\nabla w{|^2}} dxdy + \frac{m}
{2}{e^{N\theta }}\int_{{\mathbb{R}^N}} {{w^2}(x,0)} dx \hfill \\
  {\text{            }} - {e^{N\theta }}\int_{{\mathbb{R}^N}} {F(w(x,0))} dx. \hfill \\
\end{gathered} \]
By Lemma~\ref{3.1.}, $({I_m} \circ \Phi) (\theta ,w) > 0$ for all $(\theta ,w)$ with $|\theta |$, ${\| w \|_{{X^{1,s}}(\mathbb{R}_ + ^{N + 1})}}$ small and $({I_m} \circ \Phi )(0,{w_0}) < 0$, i.e. ${I_m} \circ \Phi $ possesses the Mountain-Pass geometry in $\mathbb{R} \times {X^{1,s}}(\mathbb{R}_ + ^{N + 1})$. The Mountain-Pass level of ${I_m} \circ \Phi $ is defined by
\begin{equation}\label{3.7}
{{\tilde c}_m}: = \mathop {\inf }\limits_{\tilde \gamma  \in {{\tilde \Gamma }_m}} \mathop {\sup }\limits_{t \in [0,1]} ({I_m} \circ \Phi )(\tilde \gamma (t)),
\end{equation}
where the set of paths is
\begin{equation}\label{3.8}
{{\tilde \Gamma }_m}: = \{ {\tilde \gamma  \in C([0,1],\mathbb{R} \times {X^{1,s}}(\mathbb{R}_ + ^{N + 1}) ):\tilde \gamma (0) = (0,0){\text{ and }}({I_m} \circ \Phi )(\tilde \gamma (1)) < 0} \}.
\end{equation}
As ${\Gamma _m} = \{ {\Phi  \circ \tilde \gamma :\tilde \gamma  \in {{\tilde \Gamma }_m}} \}$, the Mountain-Pass levels of ${I_m}$ and ${I_m} \circ \Phi $ coincide, i.e. ${c_m} = {{\tilde c}_m}$.

By the General Minimax principle (Theorem~2.8 of \cite{w1}), there exists a sequence $\{ ({\theta _n},{v_n})\} _{n = 1}^\infty $ in $\mathbb{R} \times {X^{1,s}}(\mathbb{R}_ + ^{N + 1})$ such that as $n \to \infty $,
\begin{equation}\label{3.9}
({I_m} \circ \Phi )({\theta _n},{v_n}) \to {c_m},
\end{equation}
\begin{equation}\label{3.10}
({I_m} \circ \Phi )'({\theta _n},{v_n}) \to 0{\text{ in (}}\mathbb{R} \times {X^{1,s}}(\mathbb{R}_ + ^{N + 1}){)^{ - 1}},
\end{equation}
\begin{equation}\label{3.11}
{\theta _n} \to 0.
\end{equation}

Indeed, set $\varepsilon  = {\varepsilon _n}: = 1/{n^2}$, $\delta  = {\delta _n}: = 1/n$ in Theorem~2.8 of \cite{w1}, \eqref{3.9}, \eqref{3.10} are direct conclusions from $(a)$, $(c)$ in Theorem~2.8 of \cite{w1}. By \eqref{3.4} and \eqref{3.5}, for $\varepsilon  = {\varepsilon _n}: = 1/{n^2}$, $\exists {\gamma _n} \in {\Gamma_m} $, such that
$\mathop {\sup }\limits_{t \in [0,1]} {I_m}({\gamma _n}(t)) \le {c_m} + 1/{n^2}$. Set ${{\tilde \gamma }_n}(t) = (0,{\gamma _n}(t))$, then
\[
\mathop {\sup }\limits_{t \in [0,1]} ({I_m} \circ \Phi) ({{\tilde \gamma }_n}(t)) = \mathop {\sup }\limits_{t \in [0,1]} {I_m}({\gamma _n}(t)) \le {c_m} + 1/{n^2}.
\]
From $(b)$ in Theorem~2.8 of \cite{w1}, $\exists ({\theta _n},{v_n}) \in \mathbb{R} \times {X^{1,s}}(\mathbb{R}_ + ^{N + 1})$ such that ${\text{dist}}(({\theta _n},{v_n}),(0,{\gamma _n}(t))) \le 2/n$, then \eqref{3.11} holds.

For every $(h,w) \in \mathbb{R} \times {X^{1,s}}(\mathbb{R}_ + ^{N + 1})$,
\begin{equation}\label{3.12}
\langle {({I_m} \circ \Phi )'({\theta _n},{v_n}),(h,w)} \rangle  = \langle {{I'_m}(\Phi ({\theta _n},{v_n})),\Phi ({\theta _n},w)} \rangle  + {P_m}(\Phi ({\theta _n},{v_n}))h.
\end{equation}
Taking $h=1$, $w=0$ in \eqref{3.12}, we have
\begin{equation}\label{3.13}
{P_m}(\Phi ({\theta _n},{v_n})) \to 0(n \to \infty) .
\end{equation}
For every $v \in {X^{1,s}}(\mathbb{R}_ + ^{N + 1})$, set $w(x,y) = v({e^{{\theta _n}}}x,{e^{{\theta _n}}}y)$, $h=0$ in \eqref{3.12}, by \eqref{3.11}, we get
\begin{equation}\label{3.14}
\langle {{I'_m}(\Phi ({\theta _n},{v_n})),v} \rangle  = o(1){\| {v({e^{{\theta _n}}}x,{e^{{\theta _n}}}y)} \|_{{X^{1,s}}(\mathbb{R}_ + ^{N + 1})}} = o(1){\| v \|_{{X^{1,s}}(\mathbb{R}_ + ^{N + 1})}}.
\end{equation}
Denote ${w_n}: = \Phi ({\theta _n},{v_n})$ in \eqref{3.9}, \eqref{3.13} and \eqref{3.14}, we get \eqref{3.6}.
\end{proof}

\begin{lemma}\label{3.3.}
Every sequence $\{ {w_n}\} _{n = 1}^\infty $ satisfying \eqref{3.6} is bounded in ${X^{1,s}}(\mathbb{R}_ + ^{N + 1})$.
\end{lemma}
\begin{proof}
By \eqref{3.6},
\begin{equation}\label{3.15}
{c_m} + o(1) = {I_m}({w_n}) - \frac{1}
{N}{P_m}({w_n}) = \frac{s}
{N}\int_{\mathbb{R}_ + ^{N + 1}} {{y^{1 - 2s}}|\nabla {w_n}{|^2}} dxdy,
\end{equation}
we get the upper bound of ${\| {{w_n}} \|_{{X^s}(\mathbb{R}_ + ^{N + 1})}}$, then by Lemma~\ref{2.2.}, we see that $\{ {w_n}(x,0)\} $ is bounded in ${L^{2_s^ * }}({\mathbb{R}^N})$. From \eqref{3.6} and \eqref{3.3}, we see that
\begin{equation}\label{3.16}
\begin{gathered}
  {\text{    }}\frac{{{k_s}(N - 2s)}}
{2}\int_{\mathbb{R}_ + ^{N + 1}} {{y^{1 - 2s}}|\nabla {w_n}{|^2}} dxdy + \frac{{mN}}
{2}\int_{{\mathbb{R}^N}} {w_n^2(x,0)} dx \hfill \\
   = N\int_{{\mathbb{R}^N}} {F({w_n}(x,0))} dx + o(1) \le \frac{{mN}}
{4}\int_{{\mathbb{R}^N}} {w_n^2(x,0)} dx + C\int_{{\mathbb{R}^N}} {w_n^{2_s^ * }(x,0)} dx + o(1), \hfill \\
\end{gathered}
\end{equation}
hence $\{ {w_n}\} $ is bounded in ${{X^{1,s}}(\mathbb{R}_ + ^{N + 1})}$.
\end{proof}
For the Mountain-Pass level $c_m$, we have the following estimate:

\begin{lemma}\label{3.4.}
If $N \ge 4s$, $2 < q < 2_s^ * $ or $2s< N < 4s$, $4s/(N - 2s) < q < 2_s^ * $, then for all $\lambda > 0$, ${c_m} < \frac{s}
{N}{({k_s}S(s,N))^{N/(2s)}}$. Moreover, if $2s< N < 4s$ and $2 < q \le 4s/(N - 2s)$, then for $\lambda  > 0$ sufficiently large, the same conclusion holds.
\end{lemma}

\begin{proof}
Let ${\phi _0} \in {C^\infty }({\mathbb{R}_ + })$ satisfying
\[
{\phi _0}(t) = 1{\text{ if }}0 \le t \le 1,{\text{ }}{\phi _0}(t) = 0{\text{ if }}t \ge 2,{\text{ }}0 \le {\phi _0}(t) \le 1{\text{ and }}|{\phi '_0}(t)| \le C
\]
and denote ${w_\delta }$ the $s$-harmonic extension of ${u_\delta }$ in Lemma~\ref{2.2.}. Denote
$${v_\delta }(x,y) = \frac{{{\phi _0}({{(|x{|^2} + {y^2})}^{1/2}}){w_\delta }(x,y)}}
{{{{\| {{\phi _0}(|x|){u_\delta }(x)} \|}_{{L^{2_ s ^*}}({\mathbb{R}^N})}}}},$$
by \cite{bcps1,dms}, we see that
\begin{equation}\label{3.17}
\| {{v_\delta }} \|_{{X^s}(\mathbb{R}_ + ^{N + 1})}^2 ={k_s} {S(s,N)} + O({\delta ^{N - 2s}}),
\end{equation}
and for any $p \in [2,2_s^ * )$,
\begin{equation}\label{3.18}
\left\| {{v_\delta }(x,0)} \right\|_{{L^p}({\mathbb{R}^N})}^p = \left\{ \begin{gathered}
  O({\delta ^{(2N - (N - 2s)p)/2}}),{\text{ if }}p > N/(N - 2s), \hfill \\
  O({\delta ^{N/2}}|\log \delta |),{\text{ if }}p = N/(N - 2s), \hfill \\
  O({\delta ^{(N - 2s)p/2}}),{\text{ if }}p < N/(N - 2s). \hfill \\
\end{gathered}  \right.
\end{equation}
By $(f_3)$,
\[\begin{gathered}
  {I_m}({v_\delta }(x/t)) \le {g_\delta }(t) \hfill \\
  : = \frac{{{k_s}}}
{2}{t^{N - 2s}}\int_{\mathbb{R}_ + ^{N + 1}} {{y^{1 - 2s}}|\nabla {v_\delta }{|^2}} dxdy - \frac{1}
{{2_s^ * }}{t^N} + \frac{m}
{2}{t^N}\int_{{\mathbb{R}^N}} {v_\delta ^2(x,0)} dx - \frac{\lambda }
{q}{t^N}\int_{{\mathbb{R}^N}} {v_\delta ^q(x,0)} dx \hfill \\
  : = {h_\delta }(t) + \frac{m}
{2}{t^N}\int_{{\mathbb{R}^N}} {v_\delta ^2(x,0)} dx - \frac{\lambda }
{q}{t^N}\int_{{\mathbb{R}^N}} {v_\delta ^q(x,0)} dx. \hfill \\
\end{gathered} \]
In view of \eqref{3.17} and \eqref{3.18}, for $\delta  > 0$ small, ${g_\delta }(t)$ has a unique critical point ${t_\delta } > 0$ which corresponds to its maximum. Therefore, we check from ${g'_\delta }({t_\delta }) = 0$ that
\[
{t_\delta } = \frac{{{{\left( {\frac{{N - 2s}}
{2}} \right)}^{1/2s}}\left\| {{v_\delta }} \right\|_{{X^s}(\mathbb{R}_ + ^{N + 1})}^{1/s}}}
{{{{\left( {\frac{{N - 2s}}
{2} - \frac{{mN}}
{2}\left\| {{v_\delta }(x,0)} \right\|_{{L^2}({\mathbb{R}^N})}^2 + \frac{{\lambda N}}
{q}\left\| {{v_\delta }(x,0)} \right\|_{{L^q}({\mathbb{R}^N})}^q} \right)}^{1/2s}}}},
\]
then we see from \eqref{3.17} and \eqref{3.18} that
\begin{equation}\label{3.19}
0 < {C_1} \le {t_\delta } \le {C_2}{\text{ for }}\delta  > 0{\text{ small.}}
\end{equation}
By \eqref{3.17} and \eqref{3.18}, we get
\begin{equation}\label{3.20}
\mathop {\max }\limits_{t \ge 0} {h_\delta }(t) = {h_\delta }({t'_\delta }) = \frac{s}
{N}{({k_s}S(s,N))^{N/2s}} + O({\delta ^{N - 2s}}),
\end{equation}
where ${t'_\delta } = \| {{v_\delta }} \|_{{X^s}(\mathbb{R}_ + ^{N + 1})}^{1/s}$ is the maximum point of ${h_\delta }(t)(t>0)$.

By \eqref{3.19} and \eqref{3.20}, we have
\begin{equation}\label{3.21}
\begin{gathered}
  {\text{    }}{c_m} \le\mathop {\sup }\limits_{t > 0} {I_m}({v_\delta }(x/t)) \le \mathop {\sup }\limits_{t > 0} {g_\delta }(t) \hfill \\
   = {h_\delta }({t_\delta }) + \frac{m}
{2}t_\delta ^N\int_{{\mathbb{R}^N}} {v_\delta ^2(x,0)} dx - \frac{\lambda }
{q}t_\delta ^N\int_{{\mathbb{R}^N}} {v_\delta ^q(x,0)} dx \hfill \\
   \le \frac{s}
{N}{({k_s}S(s,N))^{N/(2s)}} + O({\delta ^{N - 2s}}) + C\| {{v_\delta }(x,0)} \|_{{L^2}({\mathbb{R}^N})}^2 - C\lambda \| {{v_\delta }(x,0)} \|_{{L^q}({\mathbb{R}^N})}^q. \hfill \\
\end{gathered}
\end{equation}
Next, we distinguish the following cases:\\
(i) If $N>4s$, then $q > 2 > N/(N - 2s)$, by \eqref{3.18} and \eqref{3.21}, we get
\[
{c_m} \le \frac{s}
{N}{({k_s}S(s,N))^{N/(2s)}} + O({\delta ^{N - 2s}}) + O({\delta ^{2s}}) - \lambda \cdot O({\delta ^{(2N - (N - 2s)q)/2}}).
\]
In view of $(2N - (N - 2s)q)/2 < 2s < (N - 2s)$, we get the conclusion for $\delta>0$ small.\\
(ii) If $N=4s$, then $q > 2 = N/(N - 2s)$, by \eqref{3.18} and \eqref{3.21}, we have
\[
{c_m} \le \frac{s}
{N}{({k_s}S(s,N))^{N/(2s)}} + O({\delta ^{2s}}(1 + |\log \delta |)) - \lambda \cdot O({\delta ^{4s - sq}}).
\]
Since $4s-sq<2s$, we get the conclusion for $\delta>0$ small.\\
(iii) If $2s < N < 4s$ and $N/(N - 2s) < q < 2_s^ * $, we see from \eqref{3.18} and \eqref{3.21} that
\[
{c_m} \le \frac{s}
{N}{({k_s}S(s,N))^{N/(2s)}} + O({\delta ^{N - 2s}}) - \lambda  \cdot O({\delta ^{(2N - (N - 2s)q)/2}}).
\]
If $4s/(N - 2s) < q < 2_s^ * $, then $(N - 2s) > (2N - (N - 2s)q)/2$, we get the conclusion for $\delta>0$ small. If $N/(N - 2s) < q \le 4s/(N - 2s)$, then $(N - 2s) \le (2N - (N - 2s)q)/2$, we choose $\lambda  = {\delta ^{ - \theta }}$ with $\theta  > (2N - (q + 2)(N - 2s))/2 > 0$, we still get the conclusion for $\delta>0$ small.\\
(iv) If $2s < N < 4s$ and $q =N/(N - 2s)$, \eqref{3.18} and \eqref{3.21} yield
\[
{c_m} \le \frac{s}
{N}{({k_s}S(s,N))^{N/(2s)}} + O({\delta ^{N - 2s}}) - \lambda  \cdot O({\delta ^{N/2}}|\log \delta |).
\]
Since $(N - 2s) < N/2$, we choose $\lambda  = {\delta ^{ - \theta }}$ with $\theta  > 2s - (N/2)$, we get the conclusion for $\delta>0$ small.\\
(v) If $2s < N < 4s$ and $2<q <N/(N - 2s)$, \eqref{3.18} and \eqref{3.21} show that
\[
{c_m} \le \frac{s}
{N}{({k_s}S(s,N))^{N/(2s)}} + O({\delta ^{N - 2s}}) - \lambda  \cdot O({\delta ^{(N - 2s)q/2}}).
\]
We choose $\lambda  = {\delta ^{ - \theta }}$ with $\theta  > (q - 2)(N - 2s)/2$, we get the conclusion for $\delta>0$ small.
\end{proof}

\begin{lemma}\label{3.5.}
There is a sequence $\{ {x_n}\} _{n = 1}^\infty  \subset {\mathbb{R}^N}$ and $R > 0$, $\beta  > 0$ such that
\[
\int_{\Gamma _R^0({x_n})} {w_n^2(x,0)} dx \ge \beta ,
\]
where $\{ {w_n}\} _{n = 1}^\infty $ is the sequence given in \eqref{3.6}.
\end{lemma}

\begin{proof}
Assuming on the contrary that the lemma does not hold, then by Lemma~2.2 of \cite{fqt}, it follows that
\[
\int_{{\mathbb{R}^N}} {|{w_n}(x,0){|^p}} dx \to 0{\text{ as }}n \to \infty {\text{ for all }}2 < p < 2_s^ * .
\]
Since $\langle {I'_m({w_n}),{w_n}} \rangle  = o(1)$ and ${I_m}({w_n}) \to {c_m}$, by $(f_1)$ and $(f_2)$, we get
\[
\| {{w_n}} \|_{{X^s}(\mathbb{R}_ + ^{N + 1})}^2 + m\| {{w_n}(x,0)} \|_{{L^2}({\mathbb{R}^N})}^2 - \| {{w_n}(x,0)} \|_{{L^{2_s^ * }}({\mathbb{R}^N})}^{2_s^ * } = o(1),
\]
\begin{equation}\label{3.22}
\frac{1}
{2}\| {{w_n}} \|_{{X^s}(\mathbb{R}_ + ^{N + 1})}^2 + \frac{m}
{2}\| {{w_n}(x,0)} \|_{{L^2}({\mathbb{R}^N})}^2 - \frac{1}
{{2_s^ * }}\| {{w_n}(x,0)} \|_{{L^{2_s^ * }}({\mathbb{R}^N})}^{2_s^ * } = {c_m} + o(1).
\end{equation}
Let $l \ge 0$ be such that
\begin{equation}\label{3.23}
\| {{w_n}} \|_{{X^s}(\mathbb{R}_ + ^{N + 1})}^2 + m\| {{w_n}(x,0)} \|_{{L^2}({\mathbb{R}^N})}^2 \to l{\text{ and }}\| {{w_n}(x,0)} \|_{{L^{2_s^ * }}({\mathbb{R}^N})}^{2_s^ * } \to l.
\end{equation}
It is trivial that $l > 0$, otherwise ${\| {{w_n}} \|_{{X^{1,s}}(\mathbb{R}_ + ^{N + 1})}} \to 0$ as $n \to \infty $ which contradicts ${c_m} > 0$. By \eqref{3.22}, we get
\begin{equation}\label{3.24}
{c_m} = \frac{s}
{N}l.
\end{equation}
By Lemma~\ref{2.2.}, we see that
\begin{equation}\label{3.25}
\| {{w_n}} \|_{{X^s}(\mathbb{R}_ + ^{N + 1})}^2 + m\| {{w_n}(x,0)} \|_{{L^2}({\mathbb{R}^N})}^2 \ge {k_s}S(s,N)\| {{w_n}(x,0)} \|_{{L^{2_s^ * }}({\mathbb{R}^N})}^2.
\end{equation}
Letting $n \to \infty $ in \eqref{3.25}, we get $l \ge {({k_s}S(s,N))^{N/2s}}$, then by \eqref{3.24}, ${c_m} \ge \frac{s}
{N}{({k_s}S(s,N))^{N/(2s)}}$, which contradicts Lemma~\ref{3.4.}.
\end{proof}

\begin{proof}[\bf Proof of  Theorem~\ref{1.1.}]
Let $\{ {w_n}\} _{n = 1}^\infty $ be the sequence given in \eqref{3.6} and denote ${\tilde w_n}(x,y) = {w_n}(x + {x_n},y)$, where  $\{{x_n}\} _{n = 1}^\infty$ is the sequence given in Lemma~\ref{3.5.}. By Lemma~\ref{3.3.} and Lemma~\ref{3.5.}, we see that, up to a subsequence, $\exists \tilde w(x,y) \in {X^{1,s}}(\mathbb{R}_ + ^{N+1})\backslash \{ 0\} $ such that
${{\tilde w}_n}(x,y) \rightharpoonup \tilde w(x,y)$ in ${X^{1,s}}(\mathbb{R}_ + ^{N+1})$, ${{\tilde w}_n}(x,0) \to \tilde w(x,0)$ in $L_{{\text{loc}}}^p({\mathbb{R}^N})(1 \le p < 2_s^ * )$, ${{\tilde w}_n}(x,0) \to \tilde w(x,0)$ a.e. in ${\mathbb{R}^N}$ and ${\tilde w}$ satisfies \eqref{3.1}. Hence
\begin{equation}\label{3.26}
\begin{gathered}
  {b_m} \le {I_m}(\tilde w) = {I_m}(\tilde w) - \frac{1}
{N}{P_m}(\tilde w) = \frac{s}
{N}\| {\tilde w} \|_{{X^s}(\mathbb{R}_ + ^{N+1})}^2 \le \mathop {\underline {\lim } }\limits_{n \to \infty } \frac{s}
{N}\| {{{\tilde w}_n}} \|_{{X^s}(\mathbb{R}_ + ^{N+1})}^2 \hfill \\
  {\text{   }} = \mathop {\underline {\lim } }\limits_{n \to \infty } \frac{s}
{N}\| {{w_n}} \|_{{X^s}(\mathbb{R}_ + ^{N+1})}^2 = \mathop {\underline {\lim } }\limits_{n \to \infty } \Bigl[ {{I_m}({w_n}) - \frac{1}
{N}{P_m}({w_n})} \Bigr] = {c_m}. \hfill \\
\end{gathered}
\end{equation}
For any $\bar w \in {X^{1,s}}(\mathbb{R}_ + ^{N+1})\backslash \{ 0\} $ a solution of \eqref{3.1}, we set the path
\[\bar \gamma (t) = \left\{ \begin{gathered}
  \bar w(x/t,y/t),{\text{ if }}t > 0, \hfill \\
  0,{\text{ if }}t = 0. \hfill \\
\end{gathered}  \right.\]
Since
\begin{equation}\label{3.27}
{I_m}(\bar \gamma (t)) = {I_m}(\bar \gamma (t)) - \frac{1}
{N}{t^N}{P_m}(\bar w) = \Bigl( {\frac{1}
{2}{t^{N - 2s}} - \frac{{N - 2s}}
{{2N}}{t^N}} \Bigr)\| {\bar w} \|_{{X^s}(\mathbb{R}_ + ^N)}^2,
\end{equation}
there exists a ${T_0}>0$ large such that ${I_m}(\bar \gamma ({T_0})) < 0$ and ${I_m}(\bar \gamma (t))$ achieve the strict global maximum at $t=1$. By the definition of $c_m$, we see that ${I_m}(\bar w) \ge {c_m}$. Since $\bar w$ is arbitrary, we see that ${b_m} \ge {c_m}$. Hence, we conclude from \eqref{3.26} that ${I_m}(\tilde w) = {c_m} = {b_m}$ and ${I'_m}(\tilde w) = 0$. Arguing as Proposition~4.1.1 of \cite{dmv}, we see that $\tilde w \in {L^\infty }({\mathbb{R}^N})$. Since $\tilde w$ is nonnegative and nontrivial and $f$ is continuous, we can apply the Harnack's inequality in Lemma 4.9 of \cite{cs1} to conclude that $\tilde w$ is positive, that is, $\tilde w$ is in fact a positive ground state solution of \eqref{3.1}, hence, $u(x): = \tilde w(x,0)$ is a positive ground state solution of \eqref{1.1}.

\end{proof}

\end{document}